\documentclass[journal,twoside,web]{ieeecolor}
\UseRawInputEncoding 
\usepackage{generic}
\usepackage{cite}
\usepackage{amsmath,amssymb,amsfonts}
\usepackage{algorithmic}
\usepackage{graphicx}
\usepackage{textcomp}
\usepackage{amsmath}
\usepackage{latexsym, amssymb}
\usepackage{amsbsy}
\usepackage{amsopn, amstext}
\usepackage{hyperref}
\usepackage{cancel, color}
\usepackage{epstopdf}

\DeclareMathOperator{\Col}{Col}
\DeclareMathOperator{\Row}{Row}

\def\cal{\mathcal}

\def\diag{diag}
\def\ra{\rightarrow}

\def\d{\delta}

\def\D{\Delta}

\def\0{{\bf 0}}

\newcommand{\R}{{\mathbb R}}
\newtheorem{thm}{\textcolor[rgb]{0.00,0.25,0.50}{Theorem}}[section]

\newtheorem{cor}[thm]{\textcolor[rgb]{0.00,0.25,0.50}{Corollary}}
\newtheorem{dfn}[thm]{\textcolor[rgb]{0.00,0.25,0.50}{Definition}}
\newtheorem{prp}[thm]{\textcolor[rgb]{0.00,0.25,0.50}{Proposition}}
\newtheorem{exa}[thm]{\textcolor[rgb]{0.00,0.25,0.50}{Example}}
\newtheorem{rem}[thm]{\textcolor[rgb]{0.00,0.25,0.50}{Remark}}

\newtheorem{ass}[thm]{\textcolor[rgb]{0.00,0.25,0.50}{Assumption}}
\def\BibTeX{{\rm B\kern-.05em{\sc i\kern-.025em b}\kern-.08em
    T\kern-.1667em\lower.7ex\hbox{E}\kern-.125emX}}
\begin{document}
\title{A Remark on Evolution Equation of Stochastic Logical Dynamic Systems}
\author{Changxi Li, Jun-e Feng, Daizhan Cheng, and Xiao Zhang
\thanks{The research was supported by National Natural Science Foundation of China under Grant 61773371, 61877036, and the Natural Science Fund of Shandong Province under grant ZR2019MF002.}
\thanks{Changxi Li and Jun-e Feng are with School of Mathematics, Shandong University, Jinan 250100, P. R. China. (Corresponding author: Changxi Li.  E-mail: lichangxi@sdu.edu.cn, fengjune@sdu.edu.cn)}
\thanks{Daizhan Cheng and Xiao Zhang are with Key Laboratory of Systems and Control, Academy of Mathematics and Systems Sciences, Chinese Academy of Sciences, Beijing 100190, P. R. China. (E-mail: dcheng@iss.ac.cn, zhangxiao41@163.com)}}

\maketitle
\begin{abstract}
Modelling  is an essential procedure in analyzing  and controlling a given logical dynamic system (LDS). It has been proved that deterministic LDS can be modeled as a linear-like system using algebraic state space representation. However, due to the inherently non-linear,  it is difficult to obtain the algebraic expression of a stochastic LDS.   This paper provides a unified framework for transition analysis of LDSs with deterministic and stochastic dynamics. First, modelling of LDS with deterministic dynamics is reviewed. Then modeling of LDS with  stochastic dynamics is considered, and non-equivalence between subsystems and global system is proposed. Next, the reason for the non-equivalence is provided. Finally, consistency condition is presented for independent model and conditional independent model.
\end{abstract}


\section{Introduction}
\label{sec:introduction}
\IEEEPARstart{W}{ith} the development of networked systems, logical dynamic system (LDS)  has attracted a lot of attentions due to its wide background in biological system, economical system, social system etc \cite{ef14}, \cite{ma92}, \cite{sma12}. Two representative application fields of LDS are Boolean networks \cite{ef14} and  networked evolutionary games  \cite{cheng2015}.
LDSs investigate the  logic evolutionary dynamics on graphs \cite{rk02}. The appeal of LDS is that it merges the interactions between nodes and the evolutionary dynamics of nodes together, which makes it a theoretically valued model \cite{cheg2011}.

As one of most fundamental problems for dynamic systems, modeling  is the first step before analyzing  and controlling a given LDS \cite{wat16}. From a mathematical point of view, LDS can be divided into two types: (\emph{i}) LDS with deterministic dynamics: the evolutionary rule for each node is deterministic; (\emph{ii}) LDS with stochastic dynamics: the evolutionary rule for each node is stochastic. According to the type of LDS, there are different modeling methods \cite{sh10}.  Recently, a efficient mathematic tool, called \emph{semi-tensor product} (STP) of matrices, is proposed to model LDS \cite{cheng2011}.

STP is a generalization of traditional matrix product, and it has been successfully  applied  to Boolean networks \cite{cheng2011}, finite games \cite{chenggame1} and finite automaton \cite{xh18}. Recently, STP has been applied to LDS. A control Lyapunov function approach to feedback stabilization of logical control networks is proposed in \cite{lih2019}.  Event-triggered control of LDS is studied in \cite{bl18,guo17}. Controllability and observability of LDSs are considered in \cite{dch09}. Optimal control of LDS is investigated using STP in \cite{ef14}.

However,  most of existing results focused on LDS with deterministic dynamics. The reason comes from two aspects: (\emph{i}) deterministic system is easily analyzed; (\emph{ii}) deterministic LDS can be modelled as a linear-like system by \emph{algebraic state space method} using STP.  From a technical point of view, there are three key points in the process of modeling deterministic LDS into its algebraic formulation: projection, swap, and descending power, which are implemented by \emph{projection operator}, \emph{swap operator}, and \emph{power-reduced  operator}, respectively. Inspired by the successful application of algebraic formulation to deterministic LDSs, some researchers attempt to apply  it to stochastic LDSs \cite{yg18}-\cite{xd17}. However, when we shift our attention to LDS with stochastic dynamics, power-reduced  operator is not applicable any more.  Therefore, it is difficult to obtain the algebraic formulation of LDS with stochastic dynamics.
The  intrinsic cause of the difference in the modeling of  deterministic LDS and stochastic LDS is that the  latter is  \emph{inherently non-linear}.

This paper  aims at providing a unified framework for modeling LDSs with deterministic and stochastic dynamics. First, modeling of LDS with  deterministic dynamics is reviewed. Then modeling of LDS with  stochastic dynamics is considered, and non-equivalence between subsystems and global system is proposed. Next, the reason for the non-equivalence is revealed. Finally, consistency condition is presented for independent model and conditional independent model.

\emph{Contributions}: The contributions of this paper are threefold: (i) The  non-equivalence between subsystems and global system for stochastic LDS is proposed. (ii) Next, the reason for the
non-equivalence is provided. we find that stochastic LDS can be modelled as a \emph{non-homogeneous Markov chain} under  independent condition and as a  \emph{homogeneous Markov chain} under conditional independent condition. (iii) Consistency condition is presented for independent model and conditional independent model.

\emph{Notations}: $\R^n$ is denoted as the Euclidean space of all real $n$-vectors. $\mathcal{M}_{m\times n}$ is the set of $m\times n$ real matrices and
$\mathcal{D}_k:=\{1,2,\cdots,k\}$.
${\bf 1}_{m}$ is an $m$-dimensional vector with identity  entries and $I_{n}$ is an $n\times n$ identity matrix. $\Row_{r}(L)$ and $\Col_{r}(L)$ represent the $r$-th row and the $r$-th column of matrix $L$, respectively. $\Col(L)$ signifies the set of columns of $L$.
Let $\D_n:=\Col(I_n)$, $\d_n^i:= \Col_i(I_n)$.
$L=[\delta_n^{i_1}~\delta_n^{i_2}~\cdots~\delta_n^{i_r}]$ is called an $n\times r$-dimensional logical matrix, which is abbreviated as $L=\d_n[i_1,i_2,\cdots,i_r]$.
Denote by $\mathcal{L}_{s\times r}$ the set of $s\times r$ logical matrices.  $\Upsilon_{n}$ is the set of all $n$-dimensional probability vector and  $\Upsilon_{p\times q}$ is the set of all $p\times q$-dimensional column stochastic matrix.

The rest of this paper is organized as follows: Section II provides some preliminaries on semi-tensor product  of matrices and LDSs. Section III considers modeling of LDS with  deterministic dynamics. Section IV considers the  non-equivalence between subsystems and global system for stochastic LDSs. Section V and Section VI investigate the the reason for the non-equivalence and consistency condition, respectively. A brief conclusion is given in Section VII.
\section{Preliminaries}

The basic mathematical tool used in this paper is STP of matrices. Please refer to \cite{che12} for  more details.

\begin{dfn} \label{da.1} \cite{che12}
Suppose  $A\in {\cal M}_{m\times n}$, $B\in {\cal M}_{p\times q}$, and $l$ be the least common multiple of $n$ and $p$.
The STP of $A$ and $B$ is defined as
\begin{align*}
A\ltimes B:= \left(A\otimes I_{l/n}\right)\left(B\otimes I_{l/p}\right)\in {\cal M}_{ml/n\times ql/p},
\end{align*}
where $\otimes$ is the Kronecker product of matrices.
\end{dfn}

STP has the pseudo commutativity, which is shown as follows.
\begin{prp}\cite{che12} STP has the following commutativity.
\begin{enumerate}
  \item[(i)] Let  $z\in \mathbb{R}^t,~A\in M_{m\times n}$, then
$$z\ltimes A=(I_t\otimes A)\ltimes z.$$
  \item[(ii)] Let  $x\in\Delta_m,~y\in\Delta_n$ and define a  matrix $W_{[m,n]}\in \mathcal{M}_{mn\times mn}$, where
      \begin{align}
      \begin{array}{ccc}
      \Col _{(i-1)n+j}(W_{[m,n]})=\d_{mn}^{i+(j-1)m},
      \\i=1,2,\cdots, m;~ j=1,2,\cdots, n.
      \end{array}
      \end{align}
$W_{[m,n]}$ is called  $(m,n)$-th dimensional swap matrix. Then
$$x\ltimes y~=~W_{[m,n]}yx.$$
\end{enumerate}
\end{prp}

\begin{prp}\label{lemma2.2}
Let  $x\in\Delta_k$ and define a  matrix, called   power-reduced matrix
$$R_k=\diag\{\d_k^1, \d_k^2, \cdots, \d_k^k\}.$$
Then
$$x\ltimes x~=~R_kx.$$
\end{prp}

Assume $i\in\mathcal{D}_k$. By identifying $i\sim\delta_k^i$, we call $\delta_k^i$ the vector form of  integer $i$. A function $f:\prod_{i=1}^n{\cal D}_{k_i}\ra {\cal D}_{k_0}$ is called a mix-valued logical function.
\begin{prp}\cite{che12}
Let  $f:\prod_{i=1}^n{\cal D}_{k_i}\ra {\cal D}_{k_0}$ be a mix-valued logical function. Then there exists a unique matrix $M_f\in M_{k_0\times k}$, such that
$$
f(x_1,\cdots,x_n)=M_f\ltimes_{i=1}^n x_i.
$$
$M_f$ is called the structure matrix of $f$, and $k=\prod\limits_{i=1}^nk_i$.
\end{prp}

\begin{dfn}
Let $A\in \mathcal{M}_{p\times n}$ and $B\in \mathcal{M}_{q\times n}$. Then the Khatri-Rao Product of $A$ and $B$ is
  \begin{align}
  \begin{array}{ccl}
  A*B=[\Col_1(A)\ltimes \Col_1(B),\cdots,\Col_n(A)\ltimes \Col_n(B)]\\
  \in \mathcal{M}_{pq\times n}.
  \end{array}
  \end{align}
\end{dfn}

An LDS consists of two aspects: (i) an undirected graph $(N,E)$ with node set $N=\{1,2,\cdots,n\}$  and edge set $E\subset N\times N$; (ii) evolutionary rules for each node.
Let $x_i(t)\in {\cal D}_{k_i}$ be the state of node $i$ at time $t>0.$ Let $N_i$ be the neighbours of node $i$. If each node updates his state at time $t+1$ according to the state $x(t)=(x_1(t),x_2(t),\cdots,x_n(t))$ at time $t$, then the evolutionary dynamics can be described as follows
\begin{align}\label{eq1.1}
\begin{array}{l}
\begin{cases}
x_1(t+1)=f_1(\{x_i(t)\}_{i\in N_1})\\
x_2(t+1)=f_2(\{x_i(t)\}_{i\in N_2})\\
~~~~~~\vdots\\
x_n(t+1)=f_n(\{x_i(t)\}_{i\in N_n})
\end{cases}
\end{array}.
\end{align}
If $f_i:\prod_{j\in N_i}{\cal D}_{k_j}\ra {\cal D}_{k_i}$, then system (\ref{eq1.1}) is called deterministic. Otherwise, if $f_i:\prod_{j\in N_i}\Upsilon_{k_j}\ra \Upsilon_{k_i}$, then system (\ref{eq1.1}) is called stochastic.

\section{Algebraic Expression of LDS with Deterministic Dynamics}

Suppose the evolutionary rule for each node is deterministic. Consider a subset nodes $U\subseteq N$. Define a  projection matrices $\Phi_U$ as follows
$$
\Phi_U:=\otimes_{j=1}^{n}\gamma_j
$$
where
$$
\gamma_j:=
\begin{cases}
I_{k_j},~j\in U\\
{\bf 1}_{k_j}^T,~j\notin U.
\end{cases}
$$

By virtue of vector expression to $x_i$ and $x$, the subsystems (\ref{eq1.1}) can be described as follows
\begin{align}\label{eq2.1}
\begin{array}{llcc}
\begin{cases}
x_1(t+1)=M_1\ltimes_{i\in N_1}x_i(t)=M_1\Phi_{N_1}x(t)\vspace{1ex}\\
x_2(t+1)=M_2\ltimes_{i\in N_2}x_i(t)=M_2\Phi_{N_2}x(t)\vspace{1ex}\\
~~\vdots\\
x_n(t+1)=M_n\ltimes_{i\in N_n}x_i(t)=M_n\Phi_{N_n}x(t),
\end{cases}
\end{array}
\end{align}
where $x_i(t)\in\Delta_{k_i},~x(t)=\ltimes_{i=1}^nx_i(t),$ and $M_i\in\mathcal{L}_{k_i\times k}$ is the  structure matrix of $f_i$ with $k=\prod\limits_{i=1}^nk_i.$
Denote by $\hat{M}_i=M_i\Phi_{N_i}.$
Therefore, the global evolutionary equation can be obtained as
\begin{align}\label{eq2.2}
x(t+1)=Mx(t),
\end{align}
where
\begin{align}\label{eqg2.2}
M=\hat{M}_1*\hat{M}_2*\cdots*\hat{M}_n.
\end{align}

According to (\ref{eqg2.2}), global evolutionary equation (\ref{eq2.2}) can be obtained from subsystems (\ref{eq2.1}). On the other hand, if the global evolutionary equation (\ref{eq2.2}) is given, can we obtain the subsystems (\ref{eq2.1})? The following result reveals that the subsystems (\ref{eq2.1}) can be obtained from global evolutionary equation (\ref{eq2.2}).

\begin{prp}\label{p3.1}
Consider LDS  (\ref{eq1.1}) with its global evolutionary equation (\ref{eq2.2}). Then $M_i$ of subsystems (\ref{eq2.1}) can be obtained as follows
$$\hat{M}_i=\Phi_iM,~~i=1,2,\cdots,n.$$
\end{prp}
\begin{proof}
 For any $x(t)\in \Delta_{k}$, we have
\begin{align}\label{eq2.3}
\begin{array}{lllcc}
\Phi_iMx(t)&=&\Phi_ix(t+1)\vspace{1ex}&\\
&=&(\otimes_{j=1}^{n}\gamma_j)(\otimes_{i=1}^{n}x_i(t+1))\vspace{1ex}&\\
&=&x_i(t+1).&\\
\end{array}
\end{align}
The last equality comes from the following property of  Kronecker product
$$(A\otimes B)(C\otimes D)=AC\otimes BD.$$

\end{proof}

Proposition \ref{p3.1} reveals that global evolutionary equation (\ref{eq2.2}) and subsystems (\ref{eq2.1}) are equivalent for deterministic LDS. Fig. \ref{Fig1} shows the relationship between global evolution and subsystem evolution in deterministic model.

\begin{figure}[!htbp]
    \centering
    \includegraphics[width = 8cm]{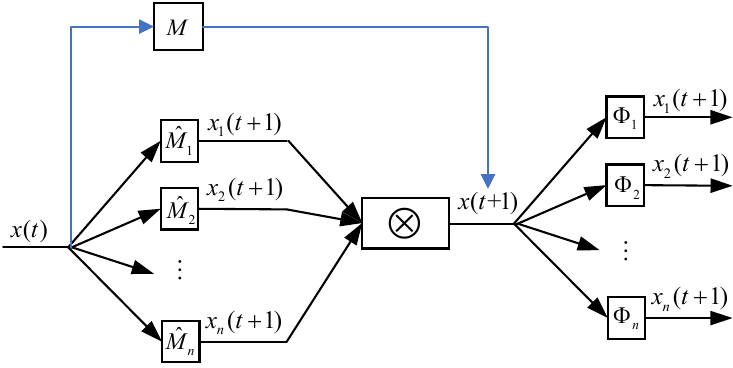}
    \caption{Evolution of state in deterministic LDS}
    \label{Fig1}  
\end{figure}

\section{Algebraic Expression of LDS with Stochastic Dynamics}

Suppose  node $\ell$ takes values from set $\mathcal{D}_{k_\ell}$ according to some probabilities  or the evolutionary  rule for each node is stochastic. Then each evolutionary  rule $f_\ell$ for node $\ell$ can be described by a transition matrix $Q_\ell=(q_{i,j}^\ell)_{k_i\times k}$ with $k=\prod\limits_{i=1}^nk_i$, where
$$q_{i,j}^\ell=\Pr(x_l(t+1)=i|x(t)=j),~\ell=1,2,\cdots,n.$$
Let $p^\ell_{j}(t)=\Pr(x_\ell(t)=j)$ be the probability of taking  $j$ of node $\ell$ at time $t.$ And the probability vector of node $\ell$ at time $t$ is denoted by
$$p_\ell(t)=[p^\ell_{1}(t),p^\ell_{2}(t),\ldots,p^\ell_{k_i}(t)]^T\in \Upsilon_{k_i}.$$
The state probability vector $p(t)$ is denoted by
$$p(t)=\big[\Pr(x(t)=1),\Pr(x(t)=2),\ldots,\Pr(x(t)=k)\big]^T.$$
Then the stochastic  evolutionary equation can be described as follows
\begin{align}\label{eq3.1.1}
\begin{array}{l}
\begin{cases}
p_1(t+1)=Q_1p_{N_1}(t)=Q_1\Phi_{N_1}p(t)\vspace{1ex}\\
p_2(t+1)=Q_2p_{N_2}(t)=Q_2\Phi_{N_2}p(t)\vspace{1ex}\\
~~\vdots\\
p_n(t+1)=Q_np_{N_n}(t)=Q_n\Phi_{N_n}p(t),
\end{cases}
\end{array}
\end{align}
where $p_{N_i}(t)\in \Upsilon_{k_{N_i}}$ is the state probability vector of players in $N_i$ and $k_{N_i}=\prod_{j\in N_i}k_j.$

Denote by $\hat{Q}_i=Q_i\Phi_{N_i}.$
Construct an   evolutionary equation, which is similar to (\ref{eq2.1}), as follows
\begin{align}\label{eq3.1}
\begin{array}{l}
\begin{cases}
p_1(t+1)=\hat{Q}_1\ltimes_{i=1}^np_i(t)\\
p_2(t+1)=\hat{Q}_2\ltimes_{i=1}^np_i(t)\\
~~\vdots\\
p_n(t+1)=\hat{Q}_n\ltimes_{i=1}^np_i(t).
\end{cases}
\end{array}
\end{align}

Construct a global evolutionary equation, which is similar to (\ref{eq2.2}), as follows
\begin{align}\label{eq3.2}
p(t+1)=Qp(t),
\end{align}
where
$$Q=\hat{Q}_1*\hat{Q}_2*\cdots*\hat{Q}_n\in \Upsilon_{k\times k}.$$

A natural question is: whether $\ltimes_{i=1}^n{p}_i(t)=p(t)$? What are the relationships between  (\ref{eq3.1.1}), (\ref{eq3.1}) and  (\ref{eq3.2})? The following example  provides some interesting results.
\begin{exa}
Consider  a two-node stochastic LDS with  transition matrix  for node $1$ and node $2$ are as follows respectively
\begin{align*}
\hat{Q}_1=\begin{bmatrix}
0.3&0.5&1&0.2\\
0.7&0.5&0&0.8
\end{bmatrix},~
\hat{Q}_2=\begin{bmatrix}
0.4&0.2&0.5&0.7\\
0.6&0.8&0.5&0.3
\end{bmatrix}
\end{align*}

According to  (\ref{eq3.1}), we have
\begin{align}\label{eq3.3}
\begin{cases}
p_1(t+1)=\hat{Q}_1\hat{p}(t),\\
p_2(t+1)=\hat{Q}_2\hat{p}(t),\\
\end{cases}
\end{align}
where $\hat{p}(t)=p_1(t)\ltimes p_2(t).$

According to  (\ref{eq3.2}), we have
\begin{align}\label{eq3.4}
\begin{array}{ccl}
p(t+1)&=&Qp(t)\vspace{1ex}\\
~&=&(\hat{Q}_1*\hat{Q}_2)p(t)\vspace{1ex}\\
~&=&\begin{bmatrix}
0.12&0.1&0.5&0.14 \\
0.18&0.4&0.5&0.06 \\
0.28&0.1&0&0.56 \\
0.42&0.4&0&0.24 \\
\end{bmatrix}p(t).
\end{array}
\end{align}

Suppose $p_1(0)=[0.4,0.6]^T,~p_2(0)=[0.5,0.5]^T,$ and
$$p(0)=p_1(0)\ltimes p_2(0)=[0.2,0.2,0.3,0.3]^T.$$
Then according to  (\ref{eq3.3}) and (\ref{eq3.4}), it is easy to calculate that
$$p_1(1)=\hat{Q}_1p(0)=[0.52,0.48]^T,$$
$$p_2(1)=\hat{Q}_2p(0)=[0.48,0.52]^T,$$
$$\hat{p}(1)=[0.2496,0.2704,0.2304,0.2496]^T,$$
$$p(1)=[0.2360,0.2840,0.2440,0.2360]^T,$$
$$\hat{p}(2)=[0.2177,0.2727,0.2262,0.2834]^T,$$
$$p(2)=[0.2118,    0.2922,    0.2266,    0.2694]^T,$$
$$\hat{p}(3)=[0.2195,    0.2650,    0.2336,    0.2819]^T,$$
$$p(3)=[0.2057,    0.2845,    0.2394,    0.2705]^T,$$
$$\vdots$$
$$\hat{p}(39)=[0.2215,    0.2665,    0.2324,    0.2796]^T,$$
$$p(39)=[0.2096,    0.2869,    0.2368,    0.2668]^T,$$
$$\hat{p}(t)=\hat{p}(39),~~p(t)=p(39),~~\forall t\geq39.$$

According to above analysis, it is easy to conclude that
$$\hat{p}(t)\neq p(t),~\forall t>0.$$

\end{exa}

Above Example implies that  subsystems (\ref{eq3.1}) and global evolutionary equation (\ref{eq3.2}) with stochastic dynamics are not equivalent. In the following section, we investigate the reasons for the non-equivalence.

\begin{rem}
The difference between (\ref{eq3.1}) and  (\ref{eq3.2}) is imperceptible, which leads to the cognitive mistake on stochastic LDS for many works. This is also the main motivation of this paper.
\end{rem}

\section{Explanation for the Non-equivalence}
Before investigating the non-equivalence, the following definitions are necessary.
\begin{dfn}
Consider three random variables $X,Y,Z$.
\begin{itemize}
\item[(i)] Random variables $X$ and $Y$ are \emph{independent} if
\begin{align*}
\begin{array}{llcc}
&\Pr(X=x,Y=y)=\Pr(X=x)\Pr(Y=y),\vspace{1ex}&\\
&~~~~~~~~~~~~~~~~~~~~~~~~~~~~~~~~~~~\forall x\in X,y\in Y.&
\end{array}
\end{align*}
\item[(ii)] Random variables $X$ and $Y$ are \emph{conditional independent} on $Z$ if
\begin{align*}
\begin{array}{lllccc}
&~~&\Pr(X=x,Y=y|Z=z)\vspace{1ex}&\\
&=&\Pr(X=x|Z=z)\Pr(Y=y|Z=z),\vspace{1ex}&\\
&~~&~~~~~~\forall x\in X,y\in Y,z\in Z.&
\end{array}
\end{align*}
\end{itemize}
\end{dfn}

\begin{rem}
It should be pointed out that neither independent  implies to conditional independent nor conditional independent implies to independent.
\end{rem}

\begin{ass}\label{as5.1} (\emph{Independence Assumption})
The random variables $x_1(t),x_2(t),\cdots,x_n(t)$ are  independent for any $t>0.$
\end{ass}

Suppose Assumption \ref{as5.1} is satisfied, then the probability of $x(t+1)=j$ is
{\small
\begin{align}\label{eq5.1}
\begin{array}{lllcc}
&~&\Pr\big(x(t+1)=j\big)\vspace{1ex}&\\
&=&\Pr\big(x_1(t+1)=j_1,\cdots,x_n(t+1)=j_n\big)\vspace{1ex}&\\
&=&\prod\limits_{i=1}^n\Pr\big(x_i(t+1)=j_i\big)\vspace{1ex}&\\
&=&\prod\limits_{i=1}^n\Pr\big(\bigcup\limits_{r=1}^k\{x(t)=r\},x_i(t+1)=j_i\big)\vspace{1ex}&\\
&=&\prod\limits_{i=1}^n\Pr\big(\bigcup\limits_{r=1}^k\{x(t)=r,x_i(t+1)=j_i\}\big)\vspace{1ex}&\\
&=&\prod\limits_{i=1}^n\Big[\sum\limits_{r=1}^k\Pr\big(x(t)=r,x_i(t+1)=j_i\big)\Big]\vspace{1ex}&\\
&=&\prod\limits_{i=1}^n\Big[\sum\limits_{r=1}^k\Pr\big(x(t)=r\big)\Pr\big(x_i(t+1)=j_i|x(t)=r\big)\Big]\vspace{1ex}&\\
&=&\prod\limits_{i=1}^n\big[\Row_{j_i}(\hat{Q}_i)p(t)\big].&\\
\end{array}
\end{align}}

According to (\ref{eq5.1}), the state probability vector can be calculated as
{\small
\begin{align}\label{eq5.2}
\begin{array}{lllccc}
\begin{bmatrix}
\Pr\big(x(t+1)=1\big)\vspace{1ex}\\
\Pr\big(x(t+1)=2\big)\vspace{1ex}\\
\vdots\\
\Pr\big(x(t+1)=k\big)
\end{bmatrix}&=&
\begin{bmatrix}
\prod\limits_{i=1}^n\Row_{1}(\hat{Q}_i)p(t)\\
\vdots\\
\prod\limits_{i=1}^n\Row_{k_i}(\hat{Q}_i)p(t)
\end{bmatrix}\vspace{1ex}&\\
&=&\ltimes_{i=1}^n[\hat{Q}_ip(t)]&
\end{array}
\end{align}}

Recall equation (\ref{eq3.1}), it is easy to know that  (\ref{eq5.2}) is equivalent to  (\ref{eq3.1}) under Assumption \ref{as5.1} (independence condition).

\begin{ass}\label{as5.2} (\emph{Conditional Independence})
The random variables $x_1(t+1),x_2(t+1),\cdots,x_n(t+1)$ are independent conditional on $x(t)$ for any $t>0.$
\end{ass}

Suppose Assumption \ref{as5.2} is satisfied, then the probability of $x(t+1)=j$ is
{\small
\begin{align}\label{eq5.3}
\begin{array}{lllcc}
&~&\Pr\big(x(t+1)=j\big)\vspace{1ex}&\\
&=&\Pr\big(\bigcup\limits_{i=1}^k\{x(t)=i\},x(t+1)=j\big)\vspace{1ex}&\\
&=&\Pr\big(\bigcup\limits_{i=1}^k\{x(t)=i,x(t+1)=j\}\big)\vspace{1ex}&\\
&=&\sum\limits_{i=1}^k\Pr\big(x(t)=i,x(t+1)=j\big)\vspace{1ex}&\\
&=&\sum\limits_{i=1}^k\Pr\big(x(t)=i\big)\Pr\big(x(t+1)=j|x(t)=i\big)\vspace{1ex}&\\
&=&\sum\limits_{i=1}^k\Big[\Pr\big(x(t)=i\big)\prod\limits_{r=1}^n\Pr\big(x(t+1)=j_r|x(t)=i\big)\Big].&\\
\end{array}
\end{align}}

According to (\ref{eq5.3}), the state probability vector can be calculated as
{\small
\begin{align}\label{eq5.4}
\begin{array}{llcc}
\begin{bmatrix}
\Pr\big(x(t+1)=1\big)\vspace{1ex}\\
\Pr\big(x(t+1)=2\big)\\
\vdots\\
\Pr\big(x(t+1)=k\big)
\end{bmatrix}
=(\hat{Q}_1*\hat{Q}_2*\cdots*\hat{Q}_n)\begin{bmatrix}
\Pr\big(x(t)=1\big)\vspace{1ex}\\
\Pr\big(x(t)=2\big)\\
\vdots\\
\Pr\big(x(t)=k\big)
\end{bmatrix}.
\end{array}
\end{align}}

Recall equation (\ref{eq3.2}), it is easy to know that  (\ref{eq5.4}) is equivalent to  (\ref{eq3.2}) under Assumption \ref{as5.2} (conditional independence condition).
\begin{figure}[!htbp]
    \centering
    \includegraphics[width = 8cm]{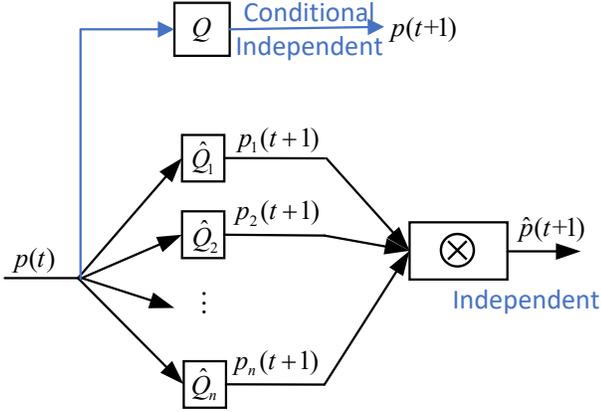}
    \caption{Evolution of state in stochastic model}
    \label{Fig2}  
\end{figure}

\begin{prp}\label{p5.1}
Consider LDS with global evolutionary equation (\ref{eq3.2}). Suppose
$$Q=\hat{Q}_1*\hat{Q}_2*\cdots*\hat{Q}_n\in \Upsilon_{k\times k}.$$
Then $\hat{Q}_i$  can be obtained as follows
$$\hat{Q}_i=\Phi_iQ,~~i=1,2,\cdots,n.$$
\end{prp}
\begin{proof}
\begin{align*}
\begin{array}{lllcc}
\Phi_iQ&=&\Phi_i[\hat{Q}_1*\hat{Q}_2*\cdots*\hat{Q}_n]\vspace{1ex}&\\
&=&\Phi_i[\otimes_{i=1}^{n}\Col_1(\hat{Q}_i),\cdots,\otimes_{i=1}^{n}\Col_{n}(\hat{Q}_i)]\vspace{1ex}&\\
&=&[\Phi_i\otimes_{i=1}^{n}\Col_1(\hat{Q}_i),\cdots,\Phi_i\otimes_{i=1}^{n}\Col_{n}(\hat{Q}_i)]\vspace{1ex}&\\
&=&[\Col_1(\hat{Q}_i),\Col_2(\hat{Q}_i),\cdots,\Col_{n}(\hat{Q}_i)]\vspace{1ex}&\\
&=&\hat{Q}_i.&\\
\end{array}
\end{align*}
\end{proof}

\begin{prp}\label{p5.1}
Consider a stochastic LDS. Then
\begin{enumerate}
  \item[(i)] (\ref{eq5.2}) is equivalent to  (\ref{eq3.1}) under independence condition.
  \item[(ii)] (\ref{eq5.4}) is equivalent to  (\ref{eq3.2}) under conditional independence condition.
  \item[(iii)] (\ref{eq3.1}) is not equivalent to  (\ref{eq3.2}).
\end{enumerate}
\end{prp}
\begin{figure}[!htbp]
    \centering
    \includegraphics[width = 5cm]{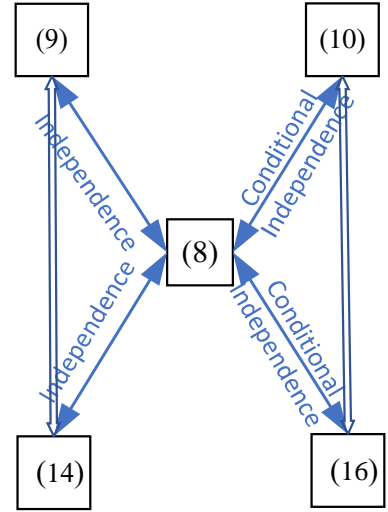}
    \caption{Relationships between different models}
    \label{Fig1}  
\end{figure}

\begin{rem}
From a mathematical point of view, (\ref{eq5.2}) is a non-homogeneous Markov chain, and (\ref{eq5.4}) is a homogeneous Markov chain.
\end{rem}

\section{Consistency Condition Exploration}

According to above analysis, we know that (\ref{eq3.1}) and (\ref{eq3.2}) are not equivalent. Under what conditions will (\ref{eq3.1}) and (\ref{eq3.2}) be  equivalent?
\begin{prp} (\emph{Consistency Condition}) \label{p5.2}
 System (\ref{eq3.1}) and (\ref{eq3.2}) are  equivalent, if and only if, for any $p\in\Upsilon_k$ the following equation is satisfied
 \begin{align}\label{eq5.5}
 HR_k^{n-1}p=H\ltimes p^{n},
 \end{align}
 where $R_k=\diag\{\delta_{k}^1,\delta_{k}^1,\cdots,\delta_{k}^k\}$ and
 $$H=\ltimes_{i=1}^n(I_{k^{i-1}}\otimes \hat{Q}_i).$$
\end{prp}
\begin{proof}
Firstly, we prove that
\begin{align}\label{eq5.1.15}
\hat{Q}_1*\hat{Q}_2*\cdots*\hat{Q}_n=H\ltimes R_k^{n-1}.
\end{align}
For any $x=\delta_k^i\in\Delta_k$
\begin{align*}
\begin{array}{lllccc}
&~~&\Col_i(\hat{Q}_1*\hat{Q}_2*\cdots*\hat{Q}_n)\vspace{1ex}&\\
&=&(\hat{Q}_1*\hat{Q}_2*\cdots*\hat{Q}_n)\ltimes x\vspace{1ex}&\\
&=&(\hat{Q}_1x)\ltimes(\hat{Q}_2x)\ltimes \cdots \ltimes (\hat{Q}_nx)\vspace{1ex}&\\
&=&\hat{Q}_1(I_k\otimes \hat{Q}_2)x^2\ltimes \cdots \ltimes (\hat{Q}_nx)\vspace{1ex}&\\
&=&Hx^n\vspace{1ex}&\\
&=&HR_k^{n-1}x.&\\
\end{array}
\end{align*}
According to the arbitrariness of $x\in \Delta_k,$ it follows that
$$\hat{Q}_1*\hat{Q}_2*\cdots*\hat{Q}_n=H\ltimes R_k^{n-1}.$$

Then we prove the following statement
\begin{align}\label{eq5.1.16}
(\hat{Q}_1p)(\hat{Q}_2p)\cdots(\hat{Q}_np)=H\ltimes p^{n},~\forall p\in\Upsilon_k.
\end{align}
For any $p\in\Upsilon_k$,
\begin{align*}
\begin{array}{lllccc}
&~~&(\hat{Q}_1p)\ltimes(\hat{Q}_2p)\ltimes \cdots \ltimes (\hat{Q}_np)\vspace{1ex}&\\
&=&\hat{Q}_1(I_k\otimes \hat{Q}_2)p^2\ltimes \cdots \ltimes (\hat{Q}_np)\vspace{1ex}&\\
&=&Hp^n.&\\
\end{array}
\end{align*}

Therefore, (\ref{eq3.1}) and (\ref{eq3.2}) are equivalent, if and only if
$$HR_k^{n-1}p=Hp^n,~~\forall p\in \Upsilon_k.$$

\end{proof}

According to Proposition \ref{p5.2}, the following corollary is obvious.
\begin{cor}
 System (\ref{eq3.1}) and (\ref{eq3.2}) are  equivalent if the following equation is satisfied
 \begin{align}\label{eq5.7}
 HR_k^{n-1}=H\ltimes p^{n-1},~~\forall p\in\Upsilon_k.
 \end{align}
\end{cor}

However, (\ref{eq5.7}) is only a sufficient condition, not a necessary condition. The following is a counterexample.
\begin{exa}
Consider  a LDS with state-strategy transition matrix as follows
\begin{align*}
\hat{Q}_1=\begin{bmatrix}
0.3&0.4&0.4&0.3\\
0.7&0.6&0.6&0.3
\end{bmatrix},~
\hat{Q}_2=\begin{bmatrix}
0.2&0.3&0.3&0.3\\
0.8&0.7&0.7&0.7
\end{bmatrix}.
\end{align*}

Let
$$p=[0,0.5,0,0.5]^T.$$

It is easy to calculate that
\begin{align*}
\begin{array}{lllcc}
HR_4\ltimes p&=&\hat{Q}_1\ltimes(I_2\otimes \hat{Q}_2)\ltimes p\vspace{1ex}&\\
&=&[0,0.35,0,0.65],\vspace{1ex}&\\
H\ltimes p^2&=& \hat{Q}_1\ltimes(I_2\otimes \hat{Q}_2)\ltimes p^2\vspace{1ex}&\\
&=&[0,0.35,0,0.65],&\\
\end{array}
\end{align*}
therefore
$$HR_4\ltimes p=H\ltimes p^2.$$

On the other hand,
\begin{align*}
\begin{array}{llcc}
HR_4-Hp=
\begin{bmatrix}
0.03&0&-0.02&0\\
0&0.02&0&-0.02\\
-0.03&0&0.02&0\\
0&-0.02&0&0.02\\
\end{bmatrix},
\end{array}
\end{align*}
therefore
$$HR_4\neq H\ltimes p.\vspace{1ex}$$
\end{exa}

\begin{prp}\label{p5.3}
 System (\ref{eq3.1}) and (\ref{eq3.2}) are  equivalent, if  there are $n-1$ matrix $\hat{Q}_{i_1},\hat{Q}_{i_2},\ldots,\hat{Q}_{i_{n-1}}$, each of which has the  same column. In other words,
 there exists  $v_i\in \Upsilon_{k_i}$ satisfying
 \begin{align}\label{eq5.8}
 \hat{Q}_{i_j}=\textbf{1}^T_{k}\otimes v_{i_j},~~j=1,2,\ldots,n-1.
 \end{align}
\end{prp}
\begin{proof}
Without loss of generality, let $\hat{Q}_1,\hat{Q}_2,\ldots,\hat{Q}_{n-1}$ be the $n-1$ matrices satisfying (\ref{eq5.8}).

 According to (\ref{eq5.1.15}),
\begin{align*}
\begin{array}{lllcc}
&~~&H\ltimes R_k^{n-1}\ltimes p\vspace{1ex}&\\
&=&(\hat{Q}_1*\hat{Q}_2*\cdots*\hat{Q}_n)p\vspace{1ex}&\\
&=&[\ltimes_{i=1}^{n}\Col_1(\hat{Q}_i),\ldots,\ltimes_{i=1}^{n}\Col_k(\hat{Q}_i)]p\vspace{1ex}&\\
&=&[\ltimes_{i=1}^{n-1}v_i\ltimes\Col_1(\hat{Q}_n),\ldots,\ltimes_{i=1}^{n-1}v_i\ltimes\Col_k(\hat{Q}_n)]p&
\end{array}
\end{align*}
 According to (\ref{eq5.1.16}),
\begin{align*}
\begin{array}{lllcc}
&~~&H\ltimes p^{n}\vspace{1ex}&\\
&=&(\hat{Q}_1p)(\hat{Q}_2p)\cdots(\hat{Q}_np)\vspace{1ex}&\\
&=&\ltimes_{i=1}^{n-1}v_i\ltimes \hat{Q}_n\ltimes p\vspace{1ex}&\\
&=&[\ltimes_{i=1}^{n-1}v_i\ltimes\Col_1(\hat{Q}_n),\ldots,\ltimes_{i=1}^{n-1}v_i\ltimes\Col_k(\hat{Q}_n)]p&\\
\end{array}
\end{align*}

Therefore,
$$HR_k^{n-1}p=H\ltimes p^{n}.$$
According to Proposition \ref{p5.2},  system (\ref{eq3.1}) and (\ref{eq3.2}) are  equivalent.

\end{proof}

\begin{exa}
Consider  a two-node stochastic LDS with  transition matrix  for node $1$ and node $2$ are as follows respectively
\begin{align*}
\hat{Q}_1=\begin{bmatrix}
0.3&0.3&0.3&0.3\\
0.7&0.7&0.7&0.7
\end{bmatrix},~
\hat{Q}_2=\begin{bmatrix}
0.2&0.6&0.1&0.4\\
0.8&0.4&0.9&0.6
\end{bmatrix}.
\end{align*}

Let
$$p=[a,b,c,1-a-b-c]^T\in\Upsilon_4.$$

It is easy to calculate that
\begin{align*}
\begin{array}{llcc}
HR_4\ltimes p&=&\hat{Q}_1\ltimes(I_2\otimes \hat{Q}_2)\ltimes p \vspace{1ex}\\
&=&\begin{bmatrix}
\frac{3}{10}a+\frac{3}{10}c \vspace{1ex}\\
\frac{3}{10}-\frac{3}{10}a-\frac{3}{10}c \vspace{1ex}\\
\frac{7}{10}a+\frac{7}{10}c \vspace{1ex}\\
\frac{7}{10}-\frac{7}{10}a-\frac{7}{10}c \vspace{1ex}\\
\end{bmatrix},\vspace{1ex}\\
H\ltimes p^2&=& \hat{Q}_1\ltimes(I_2\otimes \hat{Q}_2)\ltimes p^2\vspace{1ex}\\
&=&HR_4\ltimes p,~~~\forall p\in \Upsilon_4.\vspace{1ex}\\
\end{array}
\end{align*}

\end{exa}

The relationships between (\ref{eq3.1.1}),  (\ref{eq3.1}), (\ref{eq3.2}), (\ref{eq5.2}) and (\ref{eq5.4}) can be described as follows.
\begin{figure}[!htbp]
    \centering
    \includegraphics[width = 5cm]{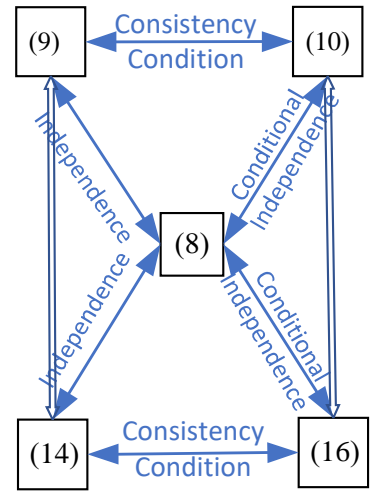}
    \caption{Relationships between different stochastic models}
    \label{Fig2}  
\end{figure}
\section{Conclusion}

As one of the fundamental problems, modeling is an essential procedure before controlling a given LDS.  This paper provided a unified framework for modeling LDSs
with deterministic and stochastic dynamics. We first reviewed the modeling of deterministic  LDS. Then non-equivalence between subsystem and global system for stochastic LDS is considered, and non-equivalence between subsystem and global system is proposed. And the reasons for the non-equivalence was provided. Finally, consistency condition was presented for independent model and conditional independent model.



\begin{thebibliography}{0}

\bibitem{ef14} E. Fornasini and M. E. Valcher,  ,``Optimal control of Boolean control networks," \emph {IEEE Trans. Autom. Control}, vol. 59, no. 5, pp. 1258--1270, 2014.

\bibitem{ma92} M. A. Nowak and R. M. May, ``Evolutionary games and spatial chaos," {\it Nature}, vol. 359, no. 6398, 826--829, 1992.

\bibitem{sma12} S. Martini,  A. Fagiolini,  L. Giarre L, et. al, ``Identification of distributed systems with logical interaction structure," {\it 51st IEEE Conference on Decision and Control}, pp. 5228--5233, 2012.

\bibitem{cheng2015} D. Cheng, F.  He, H.  Qi, and T.  Xu,``Modeling, analysis and control of networked evolutionary games,''
\emph {IEEE Trans. Autom. Control}, vol. 60, no. 9, pp. 2402--2415, 2015.


\bibitem{rk02} R. Kumar,  V. K. Garg, and S. I. Marcus, ``Finite buffer realization of input-output discrete event systems," \emph {IEEE Trans. Autom. Control}, vol. 40, no. 6, pp. 1042--1053, 2002.

\bibitem{cheg2011} Z. Yin,  Z. Li,  and D. Cheng. ``Optimal control of logical control networks," \emph {IEEE Trans. Autom. Control}, vol. 56, no. 8, pp. 1766--1776, 2011.

\bibitem{wat16}  W. Abou-Jaoud\'{e},  T. Pauline T ,  P. T. Monteiro, et. al, ``Logical modeling and dynamical analysis of cellular networks," \emph { Frontiers in Genetics}, vol. 7, no. 86, pp. 1--20, 2016.

\bibitem{sh10} S. Ilya, R. Edward, K. Seungchan, W. Zhang, ``Probabilistic Boolean networks: a rule-based uncertainty model for gene regulatory networks," {\it Bioinformatics}, vol. 18, no. 2, pp. 261--274. 2002.

\bibitem{cheng2011} D.  Cheng, H.  Qi, and Z.  Li, \emph{Analysis and Control of Boolean Networks: A Semi-tensor Product Approach}, Springer, London, 2011.

\bibitem{lih2019} H. Li, X. Ding, A control Lyapunov function approach to feedback stabilization of logical control networks. {\it SIAM Journal on Control and Optimization}, 2019, 57(2): 810--831.

\bibitem{xh18} X. Han, Z. Chen, Z. Liu, and Q. Zhang,  The detection and stabilisation of limit cycle for deterministic finite automata, {\it International Journal of Control}, vol. 91, no. 4, 874--886, 2018.

\bibitem{guo17} P. Guo, H. Zhang, F. E. Alsaadi, et al, Semi-tensor product method to a class of event-triggered control for finite evolutionary networked games, {\it IET Control
Theory \& Applications}, vol. 11, no. 13, 2140--2145, 2017.

\bibitem{chenggame1} D. Cheng, T. Xu, and H. Qi, ``Evolutionarily stable strategy of networked evolutionary games,'' \emph {	IEEE Trans. Neural Netw. Learn. Syst.},
  vol. 25, no. 7, pp. 1335--1345, 2014.

\bibitem{bl18} B. Li, L. Yang,  K. I. Kou, et. al, ``Event-triggered control for the disturbance decoupling problem of Boolean control networks," \emph {IEEE Trans. Autom. Control},
 vol. 48, no. 9, pp. 2764--2769, 2018.

\bibitem{dch09} D. Cheng and  H. Qi. ``Controllability and observability of Boolean control networks," {\it Automatica}, vol. 45, no. 7, pp. 1659--1667, 2009.

\bibitem{yg18} Y. Guo Y,  R. Zhou,  Y. Wu, et al. ``Stability and set stability in distribution of probabilistic Boolean networks," {\it IEEE Trans. Autom. Control}, vol. 64, no. 2, pp. 736--741, 2018.

\bibitem{bc22} B. Chen, J. Cao,  G. Lu, and L. Rutkowski, ``Stabilization of Markovian jump Boolean control networks via event-triggered control," {\it IEEE Trans. Autom. Control}, DOI: 10.1109/TAC.2022.3151237, 2022.

\bibitem{sz20} S. Zhu,  J. Lu,  T. Huang, et al. ``Output robustness of probabilistic Boolean control networks with respect to one-bit perturbation," {\it IEEE Transactions on Control of Network Systems}, vol. 7, no. 4, pp. 1769--1777, 2020.

\bibitem{xd17} X. Ding, H. Li, Q. Yang, Zhou Y, and F. E. Alsaadi, Stochastic stability and stabilization of $n$-person random evolutionary Boolean games, {\it Applied Mathematics and Computation}, vol. 306, 1¨C12, 2017.

\bibitem{che12} D. Cheng, H. Qi, Y. Zhao, {\it An Introduction to Semi-tensor Product of Matrices and Its Applications}, World Scientific, Singapo, 2012.

%
%
%
%
%
%
%
%
%
%
%
%
%








%
%
%
%
%
%
%
%
%
%
%
%
%
%
%
%
%
%
%
%
%
%
%




\end{thebibliography}
\end{document}